\documentclass[12pt]{amsart}
\usepackage{amsmath,amsfonts,epsfig,amssymb,graphics}

\newtheorem{thm}{Theorem}[section]
\newtheorem{cor}[thm]{Corollary}
\newtheorem{prop}[thm]{Proposition}

\newtheorem{rem}[thm]{Remark}

\newtheorem{exa}[thm]{Example}

\newcommand{\firstfunc}{\lambda_1(g) \text{Vol}(g)}
\newcommand{\kthfunc}{\lambda_k(g) \text{Vol}(g)}
\newcommand{\firstfunchigh}{\lambda_1(g) \text{Vol}(g)^{2/n}}
\newcommand{\firstfuncinv}{\lambda_1^{\text{inv}}(g) \text{Vol}(g)}
\newcommand{\kthfuncinv}{\lambda_k^{\text{inv}}(g) \text{Vol}(g)}
\newcommand{\invevk}{\lambda_k^{\text{inv}}(g)} 
 
\newcommand{\Gevk}{\lambda_k^{G}(g)} 
\newcommand{\firstfuncG}{\lambda_1^G(g) \text{Vol}(g)^{2/n}}
\newcommand{\R}{\mathbb R}

\begin{document}

\title[Extremal $G$-invariant eigenvalues]{Extremal $G$-invariant eigenvalues of the Laplacian of $G$-invariant metrics}
\author{Bruno Colbois}
\address{Universit\'e de Neuch\^atel, Institut de Math\'ematiques, Rue Emile-Argand 11, Case postale 158, 2009 Neuch\^atel
Switzerland}
\email{bruno.colbois@unine.ch}
\author{Emily B. Dryden}
\address{Departamento de Matem\'{a}tica, Instituto Superior T\'{e}cnico, Av. Rovisco Pais, 1049-001 Lisboa, Portugal}
\email{dryden@math.ist.utl.pt}
\author{Ahmad El Soufi}
\address{Laboratoire de Math\'{e}matiques et Physique Th\'{e}orique, UMR-CNRS 6083, Universit\'{e} de Tours, Parc de Grandmont, 37200 Tours, France}
\email{elsoufi@univ-tours.fr}
\date{}

\begin{abstract}
The study of extremal properties of the spectrum often involves restricting the metrics under consideration.  Motivated by the work of Abreu and Freitas in the case of the sphere $S^2$ endowed with $S^1$-invariant metrics, we consider the subsequence $\lambda_k^G$ of the spectrum of a Riemannian manifold $M$ which corresponds to metrics and functions invariant under the action of a compact Lie group $G$.  If $G$ has dimension at least 1, we show that the functional $\lambda_k^G$ admits no extremal metric under volume-preserving $G$-invariant deformations.  If, moreover, $M$ has dimension at least three, then the functional $\lambda_k^G$ is unbounded when restricted to any conformal class of $G$-invariant metrics of fixed volume.  As a special case of this, we can consider the standard $O(n)$-action on $S^n$; however, if we also require the metric to be induced by an embedding of $S^n$ in $\mathbb{R}^{n+1}$, we get an optimal upper bound on $\lambda_k^G$.  
\end{abstract}

\maketitle


\vspace{0.4cm}
Mathematics Subject Classification (2000): 58J50, 58E11, 35P15

Keywords: Laplacian, eigenvalue, invariant, extremal metric, upper bound


\section{Introduction}

\subsection{Historical background and motivation}

Let $M$ be a compact, connected manifold of dimension $n \geq 2$.  To every Riemannian metric $g$ on $M$, we can associate the Laplace-Beltrami operator $\Delta_g$ and its spectrum
$$
\text{Spec}(g) = \{ 0 = \lambda_0(g) < \lambda_1(g) \leq \lambda_2(g) \leq \cdots \leq \lambda_k(g) \leq \cdots \}.
$$
Consider the $k$th eigenvalue as a functional
$$
g \mapsto \lambda_k (g)
$$
on the space of Riemannian metrics of fixed volume on $M$; alternatively, consider the normalized functional 
$$
g \mapsto \lambda_k (g) \text{Vol}(g)^{2/n}
$$
on the space of all Riemannian metrics on $M$.  Of course, these two functionals have the same extremal metrics.

Historically, J. Hersch \cite{H} is the first to have obtained a result on these functionals; he showed that for any Riemannian metric $g$ on the sphere $S^2$, we have the relation
$$
\firstfunc \leq 8 \pi
$$
with $8 \pi = \lambda_1(g_{can})\text{Vol}(g_{can})$, where $g_{can}$ denotes the constant curvature metric on the sphere.  Moreover, the case of equality characterizes the constant curvature metric.

Following this, P. Yang and S.-T. Yau \cite{YY} showed (see also \cite{EI1}) that for an orientable surface of genus $\gamma$, we have 
$$
\firstfunc \leq 8 \pi \left[\frac{\gamma + 3}{2}\right],
$$
where $[ \cdot ]$ denotes the floor function.  This result was generalized to the nonorientable case by P. Li and Yau \cite{LY}, then to all $k$ by N. Korevaar in \cite{K}: for a compact orientable surface $M$ of genus $\gamma$, there exists a universal constant $C>0$ such that for every integer $k \geq 1$ and every Riemannian metric $g$ on $M$, we have
$$
\kthfunc \leq C(\gamma + 1) k.
$$
In general, the bound is not optimal, and finding such optimal bounds is very difficult (see \cite {N} for the case of the torus and \cite{JNP, EGJ} for that of the Klein bottle).

In dimensions three and higher, the situation differs significantly from the surface case.  
In fact, the first author and J. Dodziuk showed in \cite{CD} that, for every compact manifold $M$ of dimension $n \geq 3$, we have 
$$
\sup \firstfunchigh  = \infty,
$$
where the supremum is taken over all Riemannian metrics on $M$.  

To study extremal properties of the spectrum, it is therefore reasonable to impose additional constraints.  For example, we can restrict to a conformal class of metrics \cite{EI, K}, to projective K\"{a}hler metrics \cite{BLY}, to metrics which preserve a symplectic or K\"{a}hler structure \cite{P}, or to metrics invariant under the action of a Lie group \cite{AF, E1}.  We will focus on the last restriction; this can be examined in the context of foliations and the basic Laplacian (e.g., \cite{Ric}), although we take a different approach.  

M. Abreu and P. Freitas \cite{AF} (see also \cite{E1}) examined the case of the sphere $S^2$ endowed with $S^1$-invariant metrics.  They considered the sequence of eigenvalues of the Laplacian on $S^1$-invariant functions, $\{ \invevk, k \in \mathbb{N} \}$.  Under such constraints they showed that the functional $\firstfuncinv$ is unbounded, but becomes bounded upon restricting to the class of metrics arising from embeddings of $S^2$ in $\mathbb{R}^3$ (i.e., surfaces of revolution diffeomorphic to $S^2$).  Moreover, the supremum of the functional $\kthfuncinv$ is attained by the union of two disks of equal area, a singular surface.

This work inspires us to ask:
\begin{itemize}
\item What remains of these results when we replace $S^1$ by a general compact Lie group $G$ and $S^2$ by a general compact differentiable $G$-manifold?
\item In this general setting, do there exist extremal metrics for the functional $\kthfuncinv$?
\end{itemize}

Let $G$ be a compact Lie group acting differentiably and effectively on a compact Riemannian manifold $(M,g)$.  We assume that the action of $G$ is not transitive on $M$ and denote by
$$
\text{Spec}_G(g) = \{ 0 = \lambda_0^G(g) < \lambda_1^G(g) \leq \lambda_2^G(g) \leq \cdots \leq \lambda_k^G(g) \leq \cdots \}
$$
the subsequence of $\text{Spec}(g)$ composed of eigenvalues of the Laplacian $\Delta_g$ acting on $G$-invariant functions on $M$.  Of course, for every $k \in \mathbb{N}$, there exists an integer $m(k,g) \geq k$ such that 
$$
\Gevk = \lambda_{m(k,g)}(g).
$$
An interesting question would also be to study the behavior of $m(k,g)$ and in particular that of $m(1,g)$, which corresponds to the first energy level in which we can find eigenstates invariant by $G$.

\subsection{Nonexistence of $G$-invariant extremal metrics}

Despite the non-differentiability of the functional $g\mapsto\lambda_k^G(g)$ with respect
to metric deformations, a natural notion of extremal (or critical) metric can be
introduced. Indeed, perturbation theory enables us to prove that, for any analytic
 $G$-invariant deformation $g_\varepsilon$ of a  $G$-invariant metric $g$ (e.g., $g_\varepsilon=g+\varepsilon h$, where $h$ is a $G$-invariant symmetric tensor), the
function $\varepsilon \mapsto \lambda_k^G(g_\varepsilon)$ always admits left and right derivatives at $\varepsilon =0$ (see \cite{EI3, EI4}). The
metric $g$ is then said to be \emph{extremal} for the functional
$\lambda_k^G$ if, for any
volume preserving $G$-invariant deformation $g_\varepsilon$ of $g$, one has
$$\frac{d}{d\varepsilon } \lambda_k^G(g_\varepsilon) \Big|_{\varepsilon=0^+}\times
\frac{d}{d\varepsilon } \lambda_k^G(g_\varepsilon)
\Big|_{\varepsilon=0^-}\le 0;$$
this means that either 
$$\lambda_k^G(g_\varepsilon)\leq\lambda_k^G(g)+o(\varepsilon)
\;\;\mbox{ as } \varepsilon\rightarrow 0,$$
or
$$\lambda_k^G(g_\varepsilon)\geq\lambda_k^G(g)+o(\varepsilon)
\;\;\mbox{ as } \varepsilon\rightarrow 0.$$
In particular, if a metric $g$ is a local minimizer or a local maximizer of $\lambda_k^G$ over the set of $G$-invariant metrics of fixed volume, then $g$ is extremal in the sense of this definition.

Notice that a metric $g$ is extremal  for $\lambda_k^G$ under volume preserving $G$-invariant deformations if and only if it is extremal for $\lambda_k^G \text{Vol}^{2/n}$ under general $G$-invariant deformations. 
 \begin{thm}\label{critic} Assume that the dimension of the Lie group $G$ is at least 1. Then, for all $k\ge1$, the functional $\lambda_k^G$ admits no extremal metric under volume preserving $G$-invariant deformations. 
 \end{thm}
 Consequently, there exist neither local minimizers nor local maximizers for the functional  $\lambda_k^G$ over the set of $G$-invariant metrics of fixed volume.

\subsection{Large $G$-invariant eigenvalues in a conformal class }

According to Theorem \ref{critic}, the supremum of $\lambda_k^G(g)$ over the set $\mathcal{R}_G(M)$ of $G$-invariant metrics
 of volume 1 is either infinite, or achieved at the ``boundary'' of $\mathcal{R}_G(M)$ by a singular configuration. 
 The following theorem tells us that when the group $G$ is of dimension at least 1, only the first alternative may occur. 
 Moreover, in contrast to the functional
 $\lambda_k$, which is bounded over each conformal class of metrics (see \cite{EI,K}), the functional $\lambda_k^G$ remains
  unbounded when restricted to any conformal class of metrics in $\mathcal{R}_G(M)$.

\begin{thm}\label{grandevp}
Let $(M,g_0)$ be a compact Riemannian manifold of dimension $n \geq 3$ and $G$ a compact Lie group of dimension at least 1 
acting effectively and nontransitively on $(M,g_0)$ by isometries.  Then
$$
\sup \{ \firstfuncG \} = \infty,
$$
where the supremum is taken over all metrics $g$ which are $G$-invariant and conformal to $g_0$.
\end{thm}

\begin{rem}\label{remsphere}
As a particular case of Theorem \ref{grandevp}, one can consider metrics on the sphere $S^n$ which are invariant under the standard $O(n)$-action fixing the north and south poles.  As we will see in Theorem \ref{sphere}, the situation changes completely when we add the constraint that the metrics are induced by embeddings of $S^n$ into $\mathbb{R}^{n+1}$.
\end{rem}

As alluded to above, a result of Korevaar \cite{K} guarantees the existence of a constant $C_n([g_0])$ depending only on the dimension $n$ and the conformal class $[g_0]$ of $g_0$, so that, for any $g$ conformal to $g_0$ and any positive integer $k$,
$$
\lambda_k(g)Vol(g)^{2/n}\le C_n([g_0])k^{2/n}.
$$ 
Together with Theorem \ref{grandevp} this gives the following
\begin{cor}
Let $(M,g_0)$ and $G$ be as in Theorem \ref{grandevp}. For any positive integer $N$, there exists a G-invariant metric $g_N$ conformal to $g_0$ such that none of the first $N$ eigenfunctions of $g_N$ is G-invariant.  
\end{cor}
In other words, it is possible to conformally deform a metric $g_0$ among G-invariant metrics so that the integer $m(1,g)$ defined by  $\lambda_1^G(g)=\lambda_{m(1,g)}$ is as large as we want. 

Notice that in Theorem \ref{grandevp}, the assumption on the dimension of $G$ is necessary. Indeed, if $G$ is a finite group acting without fixed points on $M$, then any $G$-invariant metric $g$ on $M$ induces a metric  $\bar g$ on the manifold $M/G$ so that $\text{Spec}_G(g)=\text{Spec}(\bar g)$. Korevaar's result then tells us that $\lambda^G_k(g)Vol(g)^{2/n}$ is uniformly bounded over any conformal class of $G$-invariant metrics on $M$. 

Nevertheless, without the restriction to a conformal class, we still have the following 

\begin{thm}\label{casdiscret} 
Let $M$ be a compact manifold of dimension at least three and 
$G<\text{Diff}(M)$ a finite group.
Then  
$$sup \{\lambda_{1}^{G}(g)Vol(g)^{2/n} \} = \infty$$
where the supremum is taken over the set of all $G$-invariant Riemannian metrics on $M$.
\end{thm}

\subsection{Hypersurfaces of revolution}

We now return to the question of large $O(n)$-invariant eigenvalues on the sphere $S^n$ (cf. Remark \ref{remsphere}), with the additional restriction that the metrics be \emph{embedded}.  That is, let $S^n$ be embedded as a hypersurface of revolution in $(\mathbb{R}^{n+1}, \psi)$, where $\psi$ is a Riemannian metric given in Fermi coordinates
$$\R\times \R^+\times S^{n-1} \rightarrow \R^{n+1};\ \ (\rho,r,q)\rightarrow(\rho,rq)$$
by 
$$\psi(\rho,r,q)=G^2(r)d\rho^2+dr^2+F^2(r)g_0,$$
where $g_0$ denotes the canonical metric on $S^{n-1}$.

\begin{exa}
\begin{enumerate}
\item Euclidean space $\R^{n+1}$ is obtained with $G(r) = 1, F(r) = r$.
\item Hyperbolic space $\mathbb{H}^{n+1}$ is obtained with $G(r)=\cosh r, F(r) = \sinh r$.  
\end{enumerate}
\end{exa}

We will assume that $F$ and $G$ are $C^\infty$ functions and that $F$ is increasing.  This last condition is a weak one, and is clearly satisfied in the standard cases given above.  We have an isometric action of the group $O(n)$ on $(\mathbb{R}^{n+1}, \psi)$ given by
$$
(A,(\rho,r,q))\rightarrow (\rho,r,Aq)
$$
where $A \in O(n)$ and $q\in S^{n-1}$. A hypersurface in $(\mathbb{R}^{n+1}, \psi)$ which is invariant under this $O(n)$ action is said to be a ``hypersurface of revolution.''

\begin{thm}\label{sphere}
Let $g$ be an $O(n)$-invariant metric of volume 1 on $S^n$ such that $(S^n,g)$ is isometrically embedded as a  hypersurface of revolution in $(\mathbb{R}^{n+1}, \psi)$.  For all $k$, 
$$
\lambda_k^{O(n)}(g) < \lambda_k^{O(n)}(D^n) Vol (D^n)^{2/n}
$$
where $D^n$ denotes the $n$-dimensional Euclidean ball endowed with the metric $dr^2 + F^2(r)g_0$, of volume $\frac{1}{2}$, and $\lambda_k^{O(n)}(D^n)$ denotes the $k$th $O(n)$-invariant eigenvalue in the Dirichlet or Neumann spectrum of $D^n$.  

Moreover, there exists a sequence $g_i$ of $O(n)$-invariant metrics on $S^{n} \subset (\R^{n+1},\psi)$ with
$$
\lambda_{k}^{O(n)}(g_i)Vol(g_i)^{2/n} \to \lambda_k^{O(n)}(D^n)Vol(D^n)^{2/n},
$$
but the value $\lambda_k^{O(n)}(D^n)Vol(D^n)^{2/n}$ is not attained by a smooth metric on $S^n$.
\end{thm}    

The assumption that the hypersurfaces of revolution considered in this last theorem are diffeomorphic to a sphere is crucial. Indeed, in the last section, we prove that the first $S^1$-invariant eigenvalue of tori of revolution of area 1 in $\R^3$ is not bounded above.


\section{Extremal G-invariant metrics: Proof of Theorem \ref{critic}} 

Let $g$ be a $G$-invariant Riemannian metric on a compact manifold $M$ and let $k$ be a positive integer. We denote by $m$ the multiplicity of 
$\lambda_k^G(g)$, i.e. the number of  $G$-invariant eigenvalues in $spec_G(g)$ that are equal to $\lambda_k^G(g)$. For any $G$-invariant analytic
deformation $g_\varepsilon$ of $g$, one can apply general perturbation theory of unbounded self-adjoint operators to the one-parameter family of operators $\Delta_{g_\varepsilon}$; this gives the existence of a family of $m$  $G$-invariant eigenfunctions $u_{1,\varepsilon},\ldots,u_{m,\varepsilon}$ associated to a family of $m$ (unordered) $G$-invariant eigenvalues
$\Lambda_{1,\varepsilon},\ldots,\Lambda_{m,\varepsilon}$ of $(M,g_\varepsilon)$, all depending analytically on $\varepsilon$ in some interval $(-\delta, \delta)$, and satisfying
\begin{itemize}
\item  $\Lambda_{1,0}=\cdots =\Lambda_{m,0}=\lambda_k^G(  g)$, 
\item  $\forall \varepsilon \in (-\delta, \delta)$, the $m$ functions $u_{1,\varepsilon},\ldots,u_{m,\varepsilon}$ are orthonormal in $L^2(M,g_\varepsilon)$. 
\end{itemize}

From this, one can easily deduce the existence of two integers $i\le m$ and $j\le m$ such that
$$ \lambda_k^G(  g_\varepsilon) =\left\{ \begin{array}{l}
\Lambda_i(\varepsilon)\ \hbox{if}\ \varepsilon\in (-\alpha, 0)\\
\\
\Lambda_j(\varepsilon)\ \hbox{if}\ \ \varepsilon \in (0,\alpha ),
\end{array}\right.$$
for some $\alpha>0$.
Hence, the function $\varepsilon\mapsto\lambda_k^G(  g_\varepsilon)$ admits left-sided and right-sided derivatives at $\varepsilon=0$ with
$$\frac{d}{d\varepsilon} \lambda_k^G(  g_\varepsilon)\Big|_ {\varepsilon=0^-}
=  \Lambda'_i (0)$$
and
$$\frac{d}{d\varepsilon} \lambda_k^G(  g_\varepsilon) \Big|_{\varepsilon=0^+}
=  \Lambda'_j (0).$$
Following \cite{B, EI3}, if $h=\frac{d}{d\varepsilon}g_\varepsilon\Big|_{\varepsilon=0}$ is the variation tensor, then, for all $i\le m$,
$$\Lambda'_i (0)= -\int_M \langle q(u_{i,0}) , h\rangle v_g,$$
with $q(u)= du\otimes du +\frac{1}{4} \Delta_g(u^2) g$.

Let us now assume that the metric $g$ is extremal for $\lambda_k^G$ under volume-preserving $G$-invariant deformations. This implies that, for any symmetric 2-tensor 
$h$ satisfying $\int_M trace_gh \;v_g=0$, the quadratic form $Q_h:=\int_M \langle q(u) , h\rangle v_g$ is indefinite on the $G$-invariant eigenspace associated with $\lambda_k^G(g)$.

Using the same arguments as in the proof of Theorem 1.1 of \cite{EI3}, one can show that this last condition implies the existence of a finite family $f_1,\dots, f_p$ of $G$-invariant eigenfunctions associated with the eigenvalue  $\lambda_k^G(g)$ such that 
$$\sum_{i=1}^p df_i\otimes df_i = g.$$
Since $G$ is of dimension at least one and $f_1,\dots, f_p$ are constant on the orbits, $\sum_{i=1}^p df_i\otimes df_i$ cannot be a Riemannian metric on $M$, which gives the contradiction. 


\section{Large $G$-invariant eigenvalues: Proofs of Theorems \ref{grandevp} and \ref{casdiscret}}

A smooth action of a compact Lie group $G$ on a smooth manifold $M$ is well understood.  
We recall a few basic facts; for more details, see \cite[Ch. I]{Au}.
In this setting the orbits of $G$ are submanifolds of $M$, and  
$M$ has a special form in a neighborhood of an orbit.  Let $x \in M$; the slice theorem \cite[Thm. 2.1.1]{Au} 
gives the existence of an equivariant diffeomorphism from an open neighborhood of the orbit of $x$, $G \cdot x$, 
to a certain open neighborhood of $G/G_x$ ($G_x$ denotes the stabilizer of $x$ under the action of $G$).  
Furthermore, $G \cdot x$ is mapped to $G/G_x$. 
Define $V_x = T_x(M)/T_x(G \cdot x)$; then our open neighborhood of $G/G_x$ is an equivariant open neighborhood 
of the zero section in $G \times_{G_x} V_x$.  

When the quotient $M/G$ is connected, the union of the principal orbits is a dense open subset in $M$ 
\cite[Prop. 2.2.4]{Au}.  Fix a point $p$ in a principal orbit and let $U$ be an open set in $V_p$. 
By shrinking $U$ if necessary, we can identify $U$ with an open set of $\mathbb{R}^d$. 
Locally, then, an open neighborhood of $G \cdot p$ in $M$ looks like $U \times G/G_p$.  Our assumptions about the
dimension of $G$ and the nontransitivity of the action of $G$ imply that the dimension $d$ of $V_p$ as well as the dimension of the principal orbit $G \cdot p$ are between $1$ and $n-1$, and this is what we will use in the proof.

\subsection{Proof of Theorem \ref{grandevp}}

We will use the local description of $M$ given above in the following proof,  
which we complete in two steps.  We begin by assuming that the metric on $M$ can be expressed 
locally as a product metric and complete the proof in this case; we then show 
how to remove our additional assumption.             

\noindent \textbf{Step 1.}
Let $p$ be a point in a principal orbit as above.  In a neighborhood $W$ of $p$, we suppose that 
$(M,g_0)$ is $G$-equivariantly isometric to the Riemannian product $U \times (G/G_p,h)$.  
Here $U$ is endowed with a Euclidean metric and $G/G_p$ with a homogeneous metric $h$.  
Without loss of generality (up to homothety), we can assume that $Vol(G/G_p,h)=1$.  

The idea is to conformally deform this inital metric $g_0$ so that it approaches the product 
$B_{\delta} \times (G/G_p, \frac{1}{\delta^2}h)$.  Here $B_{\delta}$ is a Euclidean ball whose volume 
tends to zero with $\delta$ in such a way that $Vol(B_{\delta} \times (G/G_p, \frac{1}{\delta^2}h)) = 1$.  
We want to show that the $G$-invariant spectrum converges to the spectrum of the Euclidean ball $B_{\delta}$ 
with Neumann boundary conditions, which will imply the existence of arbitrarily large invariant first eigenvalues.    
We begin by constructing the appropriate metric on $U \times G/G_p$, and subsequently extend it to all of $M$.

Fix a point $q \in U$.  For $\delta >0$ sufficiently small, there exists $\rho = \rho (\delta)$ such that 
the Euclidean ball $B(q,\rho)$ of center $q$ and radius $\rho$ is contained in $U$ and has volume $\delta^n$.  Let $W_\delta\subset W$ be the the open set of $M$ which corresponds to $B(q,\rho) \times (G/G_p, h)$ through the local isometry above.
Thus $Vol(W_\delta,g_0) = \delta^n$ and, then,
$$Vol(W_\delta,\frac{1}{\delta^2}g_0) = 1.$$
For fixed $\delta$, we multiply the metric $\frac{1}{\delta^2}g_0$ by $\epsilon$ on the complement of 
$W_\delta$ in $M$, denoting the resulting (piecewise defined) metric on $M$ 
by $g_{\delta, \epsilon}$.  
We let $\epsilon$ tend to zero; the spectrum of $(M,g_{\delta,\epsilon})$ converges to that of 
the manifold with boundary $(W_\delta,\frac{1}{\delta^2}g_0)$ with Neumann 
boundary conditions \cite[Thm. III.1]{CV}.  The metric $g_{\delta,\epsilon}$ is of course not smooth, 
but we can make it smooth in a conformal fashion as in the proof of \cite[Thm. III.1]{CV}. 

These conformal transformations do not affect the fact that $G$ acts isometrically on $(M,g_{\delta,\epsilon})$.
In particular, we have a decomposition of the eigenfunctions of $M$ into the $G$-invariant functions 
and their orthogonal complement with respect to the $H^1$-norm.  This decomposition passes to the limit in the
theorem of Y. Colin de Verdi\`{e}re \cite{CV}.  Thus $\lambda_1^G (M,g_{\delta,\epsilon})$ converges to the first 
nonzero $G$-invariant eigenvalue of the Neumann problem on $(W_\delta,\frac{1}{\delta^2}g_0)$,
i.e., to the first eigenvalue of the $d$-dimensional ball $\left(B(q,\rho),\frac{1}{\delta^2}g_{euc}\right)$, where $g_{euc}$ is the Euclidean metric, whose volume $\delta^{n-d}$ tends to zero with $\delta$ (recall that $n>d$). 
This eigenvalue can be made arbitrarily large, concluding the argument in the case of an initial 
local product metric. 

\noindent \textbf{Step 2.}
A classical argument allows us to remove the hypothesis that $g_0$ can be expressed locally as a product.  
Using the setup given at the beginning of this section, we have an equivariant diffeomorphism from an 
open neighborhood $W$ of a point $p$ in a principal orbit of $M$ to $U \times G/G_p$.  Endowed with the 
restriction of $g_0$, $W$ is thus quasi-isometric to the Riemannian product $U \times (G/G_p,h)$, where $U$ 
is endowed with a Euclidean metric and $G/G_p$ with a homogeneous metric $h$.  The ratio of quasi-isometry 
may be very bad, but this does not matter.   We know (see \cite[p. 343]{CE}) that if two metrics are quasi-isometric, the same conformal transformation applied to both preserves the ratio of quasi-isometry; furthermore, the spectra of two quasi-isometric metrics are controlled by the ratio of quasi-isometry.  Hence for a fixed ratio of quasi-isometry, the conformal transformations in step 1 allow $\lambda_1^G$ to become arbitrarily large.

\subsection{Proof of Theorem \ref{casdiscret}}

We apply the construction of the first author and Dodziuk \cite{CD}, making it equivariant with respect to the action of $G$.  That is, let $p$ be a point in $M/G$ that is not singular.  Following \cite{CD}, we glue a sphere with a large first eigenvalue in a neighborhood of $p$ and construct an associated family of metrics.  We then lift the family of metrics to $M$ and denote it by $g_{\epsilon}$; the spectrum of $(M,g_{\epsilon})$ converges to the spectrum of a disjoint union of $|G|$ spheres.  The multiplicity of 0 is thus $|G|$.  However, the eigenvalues associated to $G$-invariant functions converge to eigenvalues associated to $G$-invariant eigenfunctions in the limit.  But the only $G$-invariant eigenfunction with eigenvalue 0 in the limit is the constant function.  Thus the first nonzero eigenvalue for $G$-invariant eigenfunctions corresponds to nonzero eigenvalues of the limit spheres and hence becomes arbitrarily large.

\begin{rem}
This suggests the following natural question: For $G$ a discrete group, can we construct $G$-invariant metrics such that the functional $\lambda_1 Vol^{2/n}$ becomes arbitrarily large?  Note that we are no longer requiring the eigenfunctions to be $G$-invariant.  The natural generalization of \cite{CD} does not extend to this case. 
\end{rem}


\section{Explicit bounds for embedded spheres: Proof of Theorem \ref{sphere}}

\subsection{Geometry of the problem}

Before going into the details of the proof of Theorem \ref{sphere}, we present the geometry of the problem.  We are concerned with spheres that are \emph{embedded} as hypersurfaces of revolution in $(\R^{n+1},\psi)$. Such a sphere of revolution  $S^n \hookrightarrow \R^{n+1}$ can be described as the result of an $O(n)$-action on a curve $c$.  A parametrization of $c$ is given by 
$$
c(t)=(\rho(t),r(t),q_0);\ \ 0<t<L,
$$
and we suppose that $c$ is parametrized by arclength, i.e. 
\begin{equation}\label{arclength}
\Vert c'(t)\Vert^2 = G^2(r(t))\rho'^2(t)+r'^2(t)=1.
\end{equation}
This implies that
$$
\vert r'(t) \vert =\sqrt{1-G^2(r(t))\rho'^2(t)}\le 1,
$$
and thus
$$ -1\le r'(t)\le 1.$$
Since our sphere is closed, we also have $r(0)=r(L)=0$, 
which implies
$$
r(t)\le t; \ \ r(t)\le L-t.
$$
Finally, the induced metric on $S^n$ is given by
$$g(t)=dt^2+ F^2(r(t))g_0,$$
and the volume form by
$$dt\wedge F^{n-1}(r(t))\omega_0,$$ 
where $\omega_0$ is the standard volume form on $S^{n-1}$.

\subsection{Proof of Theorem \ref{sphere}}

We introduce a new parametrization 
\begin{equation}
\alpha: [0,1] \rightarrow [0,L];\ s=\alpha(t)
\end{equation}
which is adapted to our problem; in particular, our parametrization is such that  
\begin{equation}\label{eqn:par}
\omega_{n-1} \int_0^{\alpha(t)} F^{n-1}(r(s))ds=t,\ 0 \le t\le 1,
\end{equation}
where $\omega_{n-1}$ denotes the volume of the unit sphere $S^{n-1}$.  The associated Riemannian metric is given by 
$$
g(r(t))= \alpha'^2(t)dt^2+F^2(r(\alpha(t)))g_0,
$$
and the volume form is, by construction,
$$
\omega(t,q)=\frac{1}{\omega_{n-1}}dt \wedge \omega_0 (q).
$$
Note that taking the derivative of both sides of (\ref{eqn:par}) gives
\begin{equation}\label{eqn:deriv}
\omega_{n-1} \alpha'(t) F^{n-1}(r(\alpha(t))) = 1.
\end{equation}

To prove our upper bounds on invariant eigenvalues, we use the Rayleigh quotient.  A function $h$ which is $O(n)$-invariant depends only on the parameter $t$, and its Rayleigh quotient is given by
\begin{eqnarray*}
R(h) & = & \frac{\int_0^1h'(t)^2 \frac{1}{\alpha'(t)^2}dt} {\int_0^1h(t)^2dt} \\
& = & \omega_{n-1}^2 \frac{\int_0^1h'(t)^2 F^{2(n-1)}(r(\alpha(t)))dt} {\int_0^1h(t)^2dt} \\
& = & \omega_{n-1}^2 \frac{\int_0^{1/2}h'(t)^2 F^{2(n-1)}(r(\alpha(t)))dt+\int_0^{1/2}h'(1-t)^2 F^{2(n-1)}(r(\alpha(1-t)))dt}{\int_0^1h(t)^2dt},
\end{eqnarray*}
where the second equality follows from (\ref{eqn:deriv}).
Recall that $F$ is an increasing function and that $r(t) \leq t$, that is, $r(\alpha(t)) \leq \alpha(t)$ and $r(\alpha (1-t)) \leq \alpha (1-t)$; these observations imply
\begin{eqnarray*}
R(h) & \leq & \omega_{n-1}^2 
\frac{\int_0^{1/2}h'(t)^2 F^{2(n-1)}( \alpha(t))dt+\int_0^{1/2}h'(1-t)^2 F^{2(n-1)}( \alpha(1-t))dt}
{\int_0^1h(t)^2dt} \\
& = & \omega_{n-1}^2 \frac{\int_0^1h'(t)^2 F^{2(n-1)}(\alpha(t))dt} {\int_0^1h(t)^2dt}.
\end{eqnarray*}

This last expression is familiar to us.  In particular, it corresponds to the case $r(t) = t$, which by equation (\ref{arclength}) implies $\rho ' (t) = 0$.  We can, without loss of generality, set $\rho(t) = 0$.  It is not hard to show that this situation corresponds exactly to the case of two balls, each of volume $\frac{1}{2}$ and with metric $dr^2 + F^2(r)g_0$, glued at their boundaries.  In particular, note that with $r(t)=t$ equation (\ref{eqn:par}) becomes 
$$
\omega_{n-1} \int_0^{\alpha(t)} F^{n-1}(s)ds=t,\ 0 \le t\le 1;
$$  
the left-hand side of this expression denotes the volume of a portion of an $n$-dimensional sphere whose boundary is a round $(n-1)$-dimensional sphere of radius $F(\alpha(t))$.  Since orthogonal projection from $\mathbb{R}^{n+1}$ to $\mathbb{R}^n$ restricted to the hypersurface decreases distance and volume, and the volume of the unit ball in $\mathbb{R}^n$ is $\frac{\omega_{n-1}}{n}$, we deduce that for $0 \leq t \leq \frac{1}{2}$, 
$$
\frac{\omega_{n-1}}{n} F^n(\alpha(t)) \leq \omega_{n-1} \int_0^{\alpha(t)} F^{n-1}(s)ds=t
$$
which implies
$$
F(\alpha(t)) \leq t^{1/n} \left(\frac{n}{\omega_{n-1}}\right)^{1/n}.
$$
The same reasoning applies to $\frac{1}{2} \leq t \leq 1$; substituting these bounds into the numerator of the Rayleigh quotient above, we see that the quotient is bounded by the expression corresponding to the case of two glued balls as claimed.  Thus, the eigenvalues (in $L^2([0,1])$) of the quadratic form  $\int_0^1h'(t)^2 \frac{1}{\alpha'(t)^2}dt$ are bounded above by those of the quadratic form obtained by gluing two balls and restricting to radial functions.  By symmetry, these latter eigenvalues are those which correspond to the radial eigenfunctions of the ball with either Dirichlet or Neumann boundary conditions.  Note that to have equality between the two quadratic forms, we must have $r(t) = t$ for all $t \in [0,1]$; this can never occur for a smooth metric.


\section{$S^1$-invariant eigenvalues of surfaces}

On closed connected surfaces there is a classification of all possible $S^1$-actions (see \cite[Chap. I, \S 3.1]{Au}).  In particular, to have a nontrivial $S^1$-action the surface must be a sphere, a torus, a projective plane, or a Klein bottle; note that these last two cannot be embedded in $\mathbb{R}^3$.  Thus to complete our study of eigenvalue behavior for $S^1$-actions on embedded surfaces we must consider the torus.  In contrast to the case of the sphere, the functional $\lambda_1^{S^1}(g) Vol(g)$ is unbounded for the torus.  

\begin{prop}
We have
$$
\sup \{\lambda_1^{S^1}(g) Vol(g)\} = \infty
$$
where the supremum is taken over all $S^1$-invariant metrics $g$ on the torus $T^2$ such that $(T^2,g)$ can be isometrically embedded as a surface of revolution in $\R^3$.
\end{prop}

\begin{rem}
This result shows that the surface being embedded in $\mathbb{R}^3$ and the metrics being invariant under an action with codimension one are not sufficient conditions for the functional $\lambda_1^{S^1}(g)Vol(g)$ to be bounded.
\end{rem}
 
\begin{proof}
We exhibit a family of $S^1$-invariant embeddings of $T^2$ in $\R^3$ with area tending to infinity and with first invariant eigenvalue uniformly bounded below by a positive constant $C$.  For this we consider the ``usual'' embeddings of the torus in $\R^3$, that is, we take a circle in the $xz$-plane with center $(R,0,0)$ and radius 1 ($R>1$) and we rotate it about the $z$-axis.  By construction, the resulting metrics are $S^1$-invariant; the volume of the resulting torus is of the order of $R$ as $R \rightarrow \infty$.  We denote this torus by $T_R$.
 
We note that $T_R$ is quasi-isometric to the product $S^1 \times S^1_R$, where $S^1$ denotes the circle of radius 1 and $S^1_R$ the circle of radius $R$.  To see this explicitly, we consider the parametrization of the torus $T_R$ given by 
$$
\psi: S^1 \times S^1 \rightarrow \R^3, \ \psi(\theta,\phi)=
((R+\cos \phi)\cos \theta, (R+\cos \phi) \sin \theta, \sin \phi),
$$ 
where the isometric action of $S^1$ is given by the translation $\theta \mapsto \theta + t$.  One can check that the associated riemannian metric can be expressed in coordinates by the diagonal matrix 
$\begin{pmatrix}
(R + \cos \phi)^2 & 0 \\
0 & 1   
\end{pmatrix}.$

We next parametrize the product $S^1 \times S^1_R$ by
$$
\alpha: S^1 \times S^1 \rightarrow \R^4;\ \alpha(\theta,\phi)=(\cos\phi,\sin \phi, R\cos \theta,R\sin \theta),
$$
where the isometric action of $S^1$ is given by the translation $\theta \mapsto \theta + t$.  One can check that the associated Riemannian metric can be expressed in coordinates by the diagonal matrix 
$\begin{pmatrix}
R^2 & 0 \\
0 & 1   
\end{pmatrix}.$

The quasi-isometry between $T_R$ and $S^1 \times S^1_R$ is obtained by associating to the point $\psi (\theta, \phi) \in T_R$, the point $\alpha (\theta, \phi) \in S^1 \times S^1_R$.  We note that this map is equivariant with respect to the isometric $S^1$-action and preserves $S^1$-invariant functions.  The ratio of quasi-isometry is given by $\frac{R^2}{(R + \cos \phi)^2}$ and this ratio is between $\frac{1}{2}$ and $\frac{3}{2}$ when $R$ is sufficiently large.

The first nonzero $S^1$-invariant eigenvalue of $S^1 \times S^1_R$ is $(2 \pi)^2$.  By a result of Dodziuk (see \cite[Prop. 3.3]{D}), $\lambda_1^{S^1}(T_R)$ is thus uniformly bounded below by a positive constant $C$ as $R \rightarrow \infty$; hence $\lambda_1^{S^1}(T_R) Vol(T_R) \rightarrow \infty$ as $R \rightarrow \infty$.
\end{proof}

In order to complete the case of $S^1$-invariant eigenvalues on the embedded torus, we note that it is not difficult to construct an $S^1$-invariant torus embedded in $\R^3$ with $\lambda_1^{S^1}(T_R) Vol(T_R) \rightarrow 0$.  To do this, we consider a ``long and thin'' ellipse given by
$$
\frac{1}{4 \epsilon ^2}(x - \epsilon)^2 + \epsilon^2 z^2
$$   
and rotate it about the $z$-axis.  Both the volume and the first $S^1$-invariant eigenvalue tend to $0$ with $\epsilon$.
  
\begin{rem}
These results on $S^1$-invariant eigenvalues on embedded tori in $\R^3$ can be extended to embedded $n$-dimensional 
products $S^1 \times S^{n-1} \subset \R^{n+1}$.
\end{rem}

\noindent \textbf{Acknowledgment.}  The second author thanks the Centre Interfacultaire Bernoulli for its support during the initial stages of this project.


\bibliographystyle{plain}
\bibliography{abcde}

\end{document}